\documentclass[12pt]{article}
\usepackage{amsmath,amssymb,amsthm} %,amsfonts,amsthm,amsxtra
\usepackage{latexcad}
\usepackage{accents}
\usepackage{bm}
\usepackage{enumitem}
\setenumerate{noitemsep,topsep=.4em,parsep=.1em,partopsep=0em} %itemsep=.1em
\def\beq#1{\begin{equation}\label{#1}}
\def\eeq{\end{equation}}
\def\beqq{\begin{eqnarray}}
\def\eeqq{\end{eqnarray}}
\def\Pr{\mathop{{\rm Pr}}\nolimits}
\def\diag{\mathop{{\rm diag}}\nolimits}
\def\l{\ell}
\def\E{\mathop{{\rm E\xy}}\nolimits}
\def\dg{\mathop{{\rm dg}}\nolimits}
\def\ldg{\bm\ell_{\dg}^\#}
\def\qdg{\bm q_{\dg}}
\def\tr{\mathop{{\rm tr}}\nolimits}
\def\cdc{,\ldots,}
\newcommand{\R}{\mathbb{R}}
\def\ll{\tilde\ell}
\def\LL{\tilde L}
\def\si{\sigma}
\def\G{\Gamma}
\def\m{\widehat m}
\def\M{\fhat M}
\def\intercal{\mathop{\scriptstyle\mathrm T}\nolimits} %\scriptscriptstyle
\def\tra{{\intercal}}
\def\x#1{} %to mark corrections; \x{??} \x{!!}
\def\Up#1{\vspace{-#1em}}
\def\xy{\hspace{.07em}}
\def\xz{\hspace{-.07em}}
\def\ms{\mathstrut}
\def\ppi{\bm\pi}
\def\proof{{\noindent\bf Proof. }}
\newtheorem{thm}{Theorem}{\bfseries}{\itshape}
\newtheorem{lem}{Lemma}{\bfseries}{\itshape}
\newtheorem{prop}{Proposition}{\bfseries}{\itshape}
\newtheorem{corol}{Corollary}{\bfseries}{\itshape}
\newtheorem{remark}{Remark}{\bfseries}{\upshape}

\DeclareMathSymbol{\widehatsym}{\mathord}{largesymbols}{"62}
\newcommand\lowerwidehatsym{%
  \text{\smash{\raisebox{-1.3ex}{%
    $\widehatsym$}}}}
\newcommand\fhat[1]{%
  \mathchoice
    {\accentset{\displaystyle\lowerwidehatsym}{#1}}
    {\accentset{\textstyle\lowerwidehatsym}{#1}}
    {\accentset{\scriptstyle\lowerwidehatsym}{#1}}
    {\accentset{\scriptscriptstyle\lowerwidehatsym}{#1}}
}
 \title{Hitting Time Quasi-metric\\ and Its Forest Representation}

\author{Pavel Chebotarev\footnote{E-mail: {\tt pavel4e@gmail.com}}\\
{\normalsize Institute of Control Sciences of the Russian Academy of Sciences}\\
{\normalsize 65 Profsoyuznaya Street, Moscow 117997, Russia}\\
Elena Deza\footnote{E-mail: {\tt Elena.Deza@gmail.com}}\\
{\normalsize Moscow State Pedagogical University}\\
{\normalsize 14 Krasnoprudnaya Street, Moscow 107140, Russia}
}

\begin{document}
\maketitle

\unitlength 1.50mm

\begin{abstract}
Let $\m_{ij}$ be the hitting (mean first passage) time from state $i$ to state $j$ in an $n$-state ergodic homogeneous Markov chain
with transition matrix~$T$. Let $\G$ be the weighted digraph %without loops
whose vertex set coincides with the set of states of the Markov chain and arc weights are equal to the corresponding transition probabilities. It holds that
$$
\m_{ij}=
q_j^{-1}\cdot
\begin{cases}
f_{ij},&\text{if }\;\; i\ne j,\\
q,     &\text{if }\;\; i=j,
\end{cases}
$$
where $f_{ij}$ is the total weight of 2-tree spanning converging forests in $\G$ that have one tree containing~$i$
and the other tree converging to $j$, $q_j$ is the total weight of spanning trees converging to~$j$ in~$\G,$
and $q=\sum_{j=1}^nq_j$ is the total weight of all spanning trees in~$\G.$
Moreover, $f_{ij}$ and $q_j$ can be calculated by an algebraic recurrent procedure.
A forest expression for Kemeny's constant is an immediate consequence of this result.
Further, we discuss the properties of the hitting time quasi-metric $m$ on the set of vertices of $\G$: $m(i,j)=\m_{ij}$, $i\neq j$, and $m(i,i)=0$.
We also consider a number of other metric structures on the set of graph vertices related to the hitting time quasi-metric $m$---along with various connections between them. The notions and relationships under study are illustrated by two examples.
\medskip

\noindent{\em Keywords:} Mean first passage time; Spanning rooted forest; Hitting time quasi-metric; Resistance metric; Commute time metric; Markov Chain Tree Theorem, Partial metric\hspace{-.3em}
\medskip

\end{abstract}

\section{Introduction}
\label{s_intr}

Let $T=[t_{ij}]\in\R^{n\times n}$ be the {\it transition matrix} of an {\it $n$-state ergodic homogeneous Markov chain}
with states $1, 2 \cdc n$. Then $T$ is an irreducible stochastic matrix.

The {\em mean first passage time} (also called the \emph{hitting time}) from state $i$ to state $j$ is %defined as follows:
\beq{defMFPT1}%+
\m_{ij}=\E(F_{ij})=\sum_{k=1}^{\infty}k\Pr(F_{ij}=k),
\eeq
where
\beq{defMFPT2}%+
F_{ij}=\min\{p>0:X_p=j\,|\,X_0=i\}
\eeq
and $X_p$ is the state of the chain at time $p\in\mathbb{N}.$
By \cite[Theorem~3.3]{Meyer75}, the matrix $\M=[\m_{ij}]\in\R^{n\times n}$ has the following representation:
\beq{Meyer3.3}%+
\M%=(I-L^\#+JL_{\dg}^\#)\Pi^{-1}
  =(I-L^\#+\bm1\ldg)\Pi^{-1},
\eeq
where $\bm1=(1\cdc 1)^\tra\in\R^n,\,$ $I=\diag(\bm1),\,$ $L^\#=[\l^\#_{ij}]_{n\times n}$ is the {\it group inverse\/} \cite{Meyer75} of $L$,
\beq{e_LaplP}
L=I-T,
\eeq
$\ldg=(\l^\#_{11}\cdc \l^\#_{nn}),$ $\Pi=\diag(\pi_1,\ldots,\pi_n)$, and $(\pi_1\cdc\pi_n)=\ppi$ is the {\it normalized left Perron vector} of $T$, i.e., the row vector in $\R^n$ satisfying
$$ %\beq{pi_vect}%-
\ppi T=\ppi \:\text{ and }\: \|\ppi\|_1=\sum_{i=1}^n\pi_i=1.
$$ %\eeq

In an entrywise form, \eqref{Meyer3.3} reads as follows (see, e.g., \cite{CatralNeumannXu05}):
\beq{Meyer3.3entry}%+
\m_{ij}=\pi_j^{-1}\cdot
\begin{cases}
(\l^\#_{jj}-\l^\#_{ij}), & \text{ if } i\ne j,\\
1,                       & \text{ if } i=j.
\end{cases}
\eeq
In the next section, we present a graph-theoretic interpretation of hitting times related to this formula.

\begin{remark}\label{r_0dia}
{\rm If one replaces $p>0$ with $p\ge0$ in the definition \eqref{defMFPT1}--\eqref{defMFPT2} of hitting time,
i.e., defines
\beq{defMFPT'}%+
m_{ij}=\E(\min\{p\ge0:X_p=j\,|\,X_0=i\}),
\eeq
then $m_{ii}=0,$ $i=1\cdc n$, and \eqref{Meyer3.3entry} and \eqref{Meyer3.3} simplify to %(\cite{CatralNeumannXu05}):
\beq{Meyer3.3entry1}%+
m_{ij}
=\frac{\l^\#_{jj}-\l^\#_{ij}}{\pi_j}
%=\pi_j^{-1}(\l^\#_{jj}-\l^\#_{ij})
,\quad i,j=1\cdc n
\eeq
and
$$ %\beq{Meyer3.3a}%-
M=[m_{ij}]_{n\times n}=(\bm1\ldg-L^\#)\Pi^{-1}.
$$ %\eeq
}
\end{remark}

\section{A forest expression for the hitting times}
\label{s_forest}

Let us say that a weighted digraph $\G$ with vertex set $V=\{1\cdc n\}$ %and without loops
{\em corresponds to the Markov chain with transition matrix\/} $T$ if
$\G$ has an arc $(i,j)$ with $i\ne j$ whenever $t_{ij}\ne0$, and the weight $w_{ij}$ of this arc is~$t_{ij}$.
Obviously, in this case the {\it Laplacian} (Kirchhoff) {\it matrix} \cite{CheAga02ap} of weighted digraph $\G,$
\beq{e_W-L}%+
\LL=\diag(W\bm1)-W,
\eeq
where $W=[w_{ij}]_{n\times n},$ coincides with $L$ of~\eqref{e_LaplP}.

Recall some graph-theoretic notation. A digraph is {\em weakly connected\/} if the corresponding undirected graph is connected. A {\em weak component\/} of a digraph $\G$ is any maximal weakly connected subdigraph of~$\G$. A {\em converging tree\/} is a weakly connected digraph in which one vertex, called the {\em root}\/, has outdegree zero and the remaining vertices have outdegree one. An {\em in-forest\/} of $\G$ is a spanning subdigraph of $\G$ all of whose weak components are converging trees (also called {\em in-arborescences}).
An in-forest is said to {\em converge to the roots\/} of its converging trees. An in-forest $F$ of a digraph $\G$ is called a {\em maximum in-forest of $\G$} if $\G$ has no in-forest with a greater number of arcs than in~$F$. The {\em in-forest connectivity of a digraph $\G$} is the number of weak components in any maximum in-forest. Obviously, every maximum in-forest of $\G$ has $n-d$ arcs, where $d$ is the in-forest connectivity of~$\G$. A~{\em submaximum in-forest of $\G$} is an in-forest of $\G$ that has $d+1$ weak components; as a consequence, it has $n-d-1$ arcs. The {\it weight of a weighted digraph} is the product of its arc weights; the weight of any digraph that has no arcs is~1. The weight of a set of digraphs is the sum of the weights of its members. In this paper, our main tool is Lemma~\ref{l_1}.

\begin{lem}[\cite{CheAga02ap}, (iii) of Proposition~15]
\label{l_1}
For any weighted digraph $\G,$ %and its Laplacian matrix~$\LL,$
it holds that
\beq{L+viaForests}%+
\LL^\#%=\si^{-1}_{n-d}\left(Q_{n-d-1}-\frac{\si_{n-d-1}}{\si_{n-d}}Q_{n-d}\right)
      =\frac{\si_{n-d-1}}{\si_{n-d}}\big(P_{n-d-1}-P_{n-d}\big),
\eeq
where $\si_k$ is the total weight of in-forests with $k$ arcs\/$,$ %(so that $\si_{n-d}$ and $\si_{n-d-1}$ are the total weights of maximum and submaximum forests of $\G$, respectively),
$P_k=Q_k/\si_k,$ and
$Q_{k}$ is the matrix whose $ij$-entry $q^{(k)}_{ij}$ $(i,j=1\cdc n)$ is the total weight of in-forests that have $k$ arcs and vertex $i$ belonging to the tree that converges to vertex~$j.$
\end{lem}

The following theorem presented in \cite{Che07mfpt} %provides
is a forest representation of the hitting times.
\begin{thm}%[\cite{Che07mfpt}]
\label{thm_MFPT}
Let $T\in\R^{n\times n}$ be the transition matrix of an $n$-state ergodic homogeneous Markov chain
with states\/ $1\cdc n$. Let $\G$ be the weighted digraph without loops whose vertices are\/ $1\cdc n$ and arc weights are
equal to the corresponding transition probabilities in~$T$.
Then the hitting time from state $i$ to state $j$ in this chain is given by
\beq{MFPT_forest}%+
\m_{ij}=
q_j^{-1}\cdot
\begin{cases}
f_{ij},&\text{if }\;\; i\ne j,\\
q,     &\text{if }\;\; i=j,
\end{cases}
\eeq
where $f_{ij}$ is the total weight of\/ $2$-tree in-forests of\/ $\G$ that have one tree containing\/ $i$ and the other tree converging to~$j,$
$q_j$ is the total weight of spanning trees converging to~$j,$ and $q={\sum_{k=1}^nq_k}=\si_{n-1}$ is the total weight of all converging trees in~$\G.$
\end{thm}

\proof
Observe that since the Markov chain under consideration is ergodic, the corresponding digraph $\G$ has a spanning converging tree. Thus, its in-forest connectivity $d$ is~1. Hence, for every $i,j=1\cdc n$, each maximum in-forest converging to $j$ is a spanning converging tree, which contains~$i$. Therefore, the $jj$- and $ij$-entries of the matrix $Q_{n-d}=Q_{n-1}$ %$=[q^{(n-1)}_{ij}]_{n\times n}$
coincide: $q^{(n-1)}_{jj}=q^{(n-1)}_{ij}=q_j$, $i, j=1\cdc n,$ where $q_j$ is the total weight of spanning trees converging to~$j.$ Thus, the differences $p^{(n-1)}_{jj}-p^{(n-1)}_{ij},$ where $[p^{(n-1)}_{ij}]_{n\times n}=P_{n-1}=P_{n-d},$ are~$0,$ as well as the differences $q^{(n-1)}_{jj}-q^{(n-1)}_{ij}.$ %$i, j=1\cdc n.$
Consequently, for any $i\ne j,$ substituting \eqref{L+viaForests} into \eqref{Meyer3.3entry} yields
\beq{eq_mijij}%+
\m_{ij}=\frac{\l^\#_{jj}-\l^\#_{ij}}{\pi_j}
       =\frac{\si_{n-2}}{\si_{n-1}\,\pi_j}\big(p^{(n-2)}_{jj}-p^{(n-2)}_{ij}\big)
       =\frac{q^{(n-2)}_{jj}-q^{(n-2)}_{ij}}{\si_{n-1}\,\pi_j}
       =\frac{f_{ij}}{\si_{n-1}\,\pi_j},
\eeq
where
\beq{fijdef}%+
f_{ij}\stackrel{\text{def}}{=}q^{(n-2)}_{jj}-q^{(n-2)}_{ij}.
\eeq
%and $[q^{(n-2)}_{ij}]_{n\times n}\!=\!Q_{n-2}\!=\!Q_{n-d-1}$.
It follows from the definition of $Q_{n-2}$ that $f_{ij}$ is the weight of the set of 2-tree in-forests of $\G$ that converge to $j$ and have $i$ and $j$ in different trees.

Furthermore, we know from the Markov Chain Tree Theorem \cite{LeightonRivest83,LeightonRivest86} obtained earlier in \cite[Lemma~7.1]{WentzellFreidlin70a} (see also \cite[Lemma~3.1]{WentzellFreidlin79e,FreidlinWentzell84} and the references in~\cite{PitmanTang16}) that
\beq{eq_MCTT}%+
\pi_j={q_j}/q,
\eeq
where %$q_j=q^{(n-1)}_{jj}$ is the total weight of trees converging to $j$ in~$\G$, so that
$q\!=\!\sum_{k=1}^nq_k\!=\!\si_{n-1}$.
Now \eqref{eq_mijij}, \eqref{eq_MCTT}, and \eqref{Meyer3.3entry} provide $\m_{ij}\!=\!\frac{f_{ij}}{q_j}$ for $i\!\ne\! j$ and $\m_{jj}\!=\!\frac q{q_j}.$\!\!\!
\qed

\begin{corol}\label{c_hi0}
For the version of hitting times introduced by \eqref{defMFPT'}$,$ in the notation of Theorem~$\ref{thm_MFPT},$ %we have
\beq{MFPT_forestA}%+
m_{ij}=
{f_{ij}}/{q_j},\quad i,j=1\cdc n
\eeq
%or in the matrix form\/$,$
and $M=\big(\bm1\qdg^{(n-2)}-Q_{n-2}\big)\big(\!\diag(q_1\cdc q_n)\big)^{-1},$ where $\qdg^{(n-2)}=\big(q_{11}^{(n-2)}\cdc q_{nn}^{(n-2)}\big).$
\end{corol}

\begin{remark}\label{r_recu}
{\rm
The values $q_j=q^{(n-1)}_{jj}$ and $f_{ij}$ that satisfy \eqref{fijdef} can be calculated by means of elementary matrix algebra, namely, by the
following recurrent procedure \cite[Proposition~4]{CheAga02ap}, which has a polynomial complexity. For $k=0,1\cdc n-2$ one has
  \beqq
  \label{req1}%+
  \si_{k+1}&=&\frac{\tr(LQ_k)}{k+1},\\
  \label{req2}%+
  Q_{k+1}&=&-LQ_k+\si_{k+1}I,
  \eeqq
where $\si_0=1$, and  $Q_0=I$.
Relations with the GTH algorithm can be found in \cite{Sonin99}.
}
\end{remark}
%Relations with the State Reduction/GTH algorithm can be found in \cite{Sonin99}.

\begin{remark}\label{r_proof}
{\rm
Theorem~\ref{thm_MFPT} can be alternatively derived from \cite[Lemma~3.3]{OlivieriScoppola96} or \cite[Lemma~3.4]{Catoni99}, both based on Lemma~3.4 in~\cite{WentzellFreidlin79e,FreidlinWentzell84}. The authors are grateful to Raphael Cerf for pointing out Ref.\:\cite{Catoni99}. Some special cases of Theorem~\ref{thm_MFPT} were obtained %independently
in \cite{Hunter16}; see also~\cite{PitmanTang16}.
}
\end{remark}

The following corollary (recently appeared as \cite[Corollary\:1.4]{PitmanTang16} and \cite[Theorem\:2.3]{KirklandZeng16}) relates to Kemeny's constant (see \cite{Hunter14}) and is an immediate consequence of Theorem~\ref{thm_MFPT}.

\begin{corol}\label{c_kemeny}
Under the conditions of Theorem~$\ref{thm_MFPT},$ Kemeny's constant is equal to $1+\dfrac{\si_{n-2}}{\si_{n-1}}.$
\end{corol}

\proof By Theorem~\ref{thm_MFPT} and \eqref{eq_MCTT}, for any $i\in V,$
\beq{e_prKem}
  \sum_{j=1}^n\pi_j\xy\m_{ij}
=                  \frac{q_i}q\!\cdot\!\frac q{q_i}
+     \sum_{j\ne i}\frac{q_j}q\!\cdot\!\frac{f_{ij}}{q_j}
= 1 + \sum_{j\ne i}\frac{f_{ij}}q
= 1 +              \frac{\si_{n-2}}{\si_{n-1}}.
\eeq
The last transition holds because the weight of any 2-tree in-forest of $\G$ is included in exactly one sum $f_{ij}$:
it is the one where $j$ is the root of the tree that does not contain~$i.$
\qed

\begin{remark}\label{r_mvsm}
{\rm
By \eqref{e_prKem} for any $i\in V,$ $\sum_{j=1}^n\pi_j\xy m_{ij} = \frac{\si_{n-2}}{\si_{n-1}}.$ Along with \eqref{Meyer3.3entry1} and Corollary~\ref{c_hi0}, this convinces that the $m_{ij}$'s lead to simpler expressions than the $\m_{ij}$'s do.
}
\end{remark}

In Sections~\ref{s_quasi} and \ref{s_simp-walk} we study metric properties of the hitting times.

\section{Hitting time quasi-metric and related metrics}
\label{s_quasi}

\subsection{Hitting time quasi-metric}
\label{ss_hitting}
A function $d\!:\! X\times X\to\R$ is a {\em quasi-metric on $X$} \cite{Hausdorff1927,Wilson1931,DezaDeza16EShort} if for all $x,y,z\in X,$
\begin{enumerate}
\item $d(x,y)\ge 0;$
\item $d(x,y)=0$ if and only if $x=y;$
\item $d(x,y)\le d(x,z)+d(z,y)$ (\xz{\it oriented triangle inequality}\xy).
\end{enumerate}

As distinct from metrics, quasi-metrics are not generally symmetric.

It was observed in Remark~\ref{r_mvsm} that the $m_{ij}$ version \eqref{defMFPT'} of hitting times leads to more elegant expressions then the $\m_{ij}$ version~\eqref{defMFPT1}.
In addition, by \cite[Proposition~9-58]{KemenySnellKnapp76}, %\cite[Theorem\:6.2.1]{KirklandNeumann12book}
$m(i,j)=m_{ij}$ is a quasi-metric on the set of states of our Markov chain (and on the set of vertices of any corresponding weighted digraph~$\G$). It is called the {\em hitting time} (or {\em mean first-passage time}) {\em quasi-metric}.

Moreover, by \cite[Proposition~9-58]{KemenySnellKnapp76} or \cite[Theorem\:6.2.1]{KirklandNeumann12book}, this quasi-metric satisfies the {\it cutpoint additivity} \cite{Che13Paris} (also called the {\it graph-geodetic property} \cite{KleinZhu98}):
$$
m(i,j)=m(i,k)+m(k,j)
$$
holds true if and only if all paths in $\G$ from $i$ to $j$ pass through~$k.$

\subsection{Commute time metric}
The {\it commute time metric} (or {\it random roundtrip time distance}) $c$ on the set of states of our Markov chain (or on $V(\G),$ where $\G$ is any corresponding weighted digraph) is defined by
\beq{e_cij}%+
c(i,j)=m(i,j)+m(j,i),\quad i,j=1\cdc n.
\eeq

The commute time $c(i,j)$ is the average number of steps that takes a random walk to reach $j$ from $i$ and return to~$i$.
$C=[c_{ij}]$ is the corresponding matrix. %has $c_{ij}=c(i,j).$
Using %\eqref{Meyer3.3entry1} and
\eqref{MFPT_forestA} we have
\begin{corol}
\label{c_commu}
Under the conditions of Theorem~$\ref{thm_MFPT},$ for all $i,j\in V,$
$
c(i,j)%=\frac{\l^\#_{jj}-\l^\#_{ij}}{\pi_j}
      %+\frac{\l^\#_{ii}-\l^\#_{ji}}{\pi_i}
      =\frac{f_{ij}}{q_j}
      +\frac{f_{ji}}{q_i}.
$
\end{corol}
Since $m(i,j)$ is a cutpoint additive quasi-metric and $c(i,j)$ is symmetric, $c(i,j)$ is a cutpoint additive metric.

\subsection{Resistance distance}
There is a strong connection between random walks in graphs and electric networks~\cite{DoyleSnell84}.
Given a \emph{connected weighted undirected\/} graph $G$, the underlying electrical network is the network obtained by replacing vertices and edges by nodes and electrical resistors, respectively. Edge weights are interpreted as conductances, so the resistances are the reciprocal weights. The {\it effective resistance\/} $\Omega(i,j)=\Omega_{ij}$ between any two nodes $i$ and $j$ is defined as the voltage that develops between $i$ and $j$ when a unit current is maintained through them (i.e., enters one and leaves the other node).

Obviously, for all nodes $i,j,k$, $\Omega(i,j)\ge0$, $\Omega(i,j)=0$ iff $i=j,$ $\Omega(i,j)=\Omega(j,i),$ and it can be shown that
$$
\Omega(i,j)+\Omega(j,k)\ge \Omega(i,k),
$$
i.e., $\Omega(\cdot,\cdot)$ is a metric \cite{Sharpe67a,GvishianiGurvich87En} called the {\it electric metric} (or the {\it %$[$effective\/$]$ 
resistance distance}~\cite{KleinRandic93}).

Let $\LL$ be the symmetric {\it Laplacian matrix\/} of $G$ defined by \eqref{e_W-L},
%\beq{e_W-L}\LL=\diag(W\bm1)-W,\eeq
where $W$ is the matrix of edge weights of~$G$. The tilde distinguishes this matrix from $L$ of~\eqref{e_LaplP}.
The resistance distance in $G$ can be represented as follows \cite{SharpeStyan67,RaoMitra71}:
\beq{e_res-inv}%+
\Omega(i,j)=\ll_{ii}^\#+\ll^\#_{jj}-\ll^\#_{ij}-\ll^\#_{ji},
\eeq
where $\LL^\#=[\ll^\#_{ij}]_{n\times n}$ is the group inverse (coinciding in this case with the Moore-Penrose generalized inverse) of~$\LL.$

Consider a forest representation of the resistance distance (\cite[Theorem~7-4]{SeshuReed61}; \cite{Shapiro87MathMag}).
\begin{corol}\label{c_res-fo}
$\Omega(i,j)=f'_{ij}/q',$ where $q'$ is the total weight of spanning trees in $G$ and $f'_{ij}$ is the total weight of $2$-tree spanning forests of $G$ having $i$ and $j$ in different trees.
\end{corol}
\proof Let $\G$ be the directed version of $G$: for every edge of $G,\xy$ $\G$ has a pair of opposite arcs carrying the weight of that edge. $G$ and $\G$ share the same~$\LL.$ In the same way as in \eqref{eq_mijij}, for $\G$ we get $\ll_{ii}^\#+\ll^\#_{jj}-\ll^\#_{ij}-\ll^\#_{ji}=(f_{ij}+f_{ji})/q.$ Returning to $G$ observe that $q=nq'$ (as any spanning tree in~$G$ corresponds to $n$ trees of the same weight converging to different vertices in $\G$) and $f_{ij}+f_{ji}=nf'_{ij}$ (as any 2-tree spanning forest in $G$ with $n_i$ vertices in the tree containing $i$ and $n_j=n-n_i$ vertices in the tree containing $j$ corresponds to $n_i$ converging forests whose weights are counted in $f_{ij}$ and $n_j$ forests whose weights are counted in $f_{ji}$). Therefore by \eqref{e_res-inv}, $\Omega(i,j)=(f_{ij}+f_{ji})/q=f'_{ij}/q'.$
\qed

There are two popular ways of attaching a Markov chain to a weighted graph~$G.$ %(di)graph with positive edge (arc) weights.
The first one is to define the transition matrix (cf.\ the first paragraph of Section\;\ref{s_forest}) by
\beq{e_L-T}%+
T_\tau=I-\tau \LL,
\eeq
where\footnote{Sometimes $\tau=(\max_i\sum_{j\ne i}w_{ij})^{-1}$ or $\tau=((n-1)\max_{i,j}w_{ij})^{-1}$ or $\tau=(n\max_{i,j}w_{ij})^{-1}$ is chosen.} $0\!<\!\tau\!\le\!(\max_i\sum_jw_{ij})^{-1},$ which guarantees the stochasticity of $T\!=\!T_\tau.$ %$T_\tau.$
Then $T$ is symmetric for any undirected~$G.$ Moreover, all transition probabilities between distinct vertices are proportional to the %corresponding
edge weights in~$G,$ as the matrix \eqref{e_LaplP} of $T$ is proportional to~$\LL.$
On the other hand, $T$ normally has a nonzero diagonal even when $G$ has no loops, which allows the corresponding Markov chain to preserve its state on adjacent steps.

The second way is to normalize each row of $W$ separately:
\beq{e_W-T}%+
T_W=(\diag(W\bm1))^{-1}W.
\eeq
Here, the symmetry of $W$ does not guarantee the symmetry of~$T=T_W,$ while the chain alters its state on each step whenever $G$ has no loops.

It is noteworthy that with either way of defining $T,$ the resistance distance for $G$ is {\em proportional\/} to the commute time metric for the Markov chain determined by~$T.$

\begin{prop}
\label{p_recoLa}
For transition matrices \eqref{e_L-T}$,$ $C=n\tau^{-1}\,\Omega,$ where $C=[c_{ij}],$ $\Omega=[\Omega_{ij}].$
%$c(i,j)=n\tau^{-1}\,\Omega(i,j)$ for all $i,j\in V(G).$
\end{prop}

\proof
Observe that by \eqref{e_LaplP}, %and \eqref{e_W-L},
$L=\tau\LL.$ Since $T$ is symmetric, $\ppi=n^{-1}\bm1^{\xz\tra}$ holds. Now comparing \eqref{Meyer3.3entry1} and \eqref{e_cij} with \eqref{e_res-inv} we have
$%\beq{e_resi-commu1}%-
c(i,j)=n\big(\l^\#_{jj}-\l^\#_{ij}+\l^\#_{ii}-\l^\#_{ji}\big)
      =n\tau^{-1}\,\Omega(i,j). %, \quad i,j=1\cdc n.
$%\eeq
\qed

\begin{corol}[\cite{ChandraRaghavanRuzzoSmolenskyTiwari89}]
\label{c_recoNo}
For the transition matrix \eqref{e_W-T}$,$ $C=\big(\sum_{k,t=1}^nw_{kt}\big)\,\Omega.$ %,\,$ $i,j\in V(G).$
%$c(i,j)=\big(\sum_{k,t=1}^nw_{kt}\big)\,\Omega(i,j),\,$ $i,j\in V(G).$
\end{corol}

\proof
By Theorem\:\ref{thm_MFPT}, $m(i,j)=\frac{f_{ij}}{q_j}$ $(i\ne j)$. Every spanning tree converging to $j$ in $\G$ has one arc weight in each row of $T,$ except for row~$j.$ Hence by \eqref{e_W-T}, $q_j=q's_jR,$ where $q'$ is the total weight of spanning trees in $G,$ $\bm s=(s_1\cdc s_n)^\tra=W\bm1,$ and $R=\big(\prod_{k=1}^ns_k\big)^{-1}.$
Every 2-tree in-forest whose weight is a term of the sum $f_{ij}$ has one arc weight in each row of $T,$ except for row~$j$ and some other row~$k.$
Hence by \eqref{e_W-T}, $f_{ij}=\sum_{k\ne j}f'_{ik,j}s_js_kR,$ where $f'_{ik,j}$ is the total weight of 2-tree spanning forests in $G$ having $i$ and $k$ in one tree and $j$ in the other tree. Therefore, $m_{ij}=\frac{f_{ij}}{q_j}=\frac{\sum_{k\ne j}f'_{ik,j}s_k}{q'}$ and $c_{ij}=m_{ij}+m_{ji}=(q')^{-1}\big(\sum_{k\ne j}f'_{ik,j}s_k+\sum_{k\ne i}f'_{jk,i}s_k\big).$ The weight of each 2-tree spanning forest of $G$ having $i$ and $j$ in different trees is a term of the first sum on the r.h.s.\ with multipliers $s_k$ for all vertices $k$ of its tree containing $i$ and enters the second sum with multipliers $s_k$ for all vertices $k$ of its tree containing~$j.$ Thus, %sum of such multipliers is $\sum_{k=1}^ns_k$ and
using Corollary\:\ref{c_res-fo} we have $c(i,j)=(q')^{-1}f'_{ij}\sum_{k=1}^ns_k=\big(\sum_{k,t=1}^nw_{kt}\big)\,\Omega(i,j).$
\qed

For additional relations between the electric metric and Markov chains, we refer to~\cite{EllensSpieksma11} and for relevant identities to~\cite{BapatSivasubramanian11}.
In \cite{YoungScardoviLeonard16}, effective resistance is generalized to directed graphs.
%\nocite{KleinRandic93,Klein97,Bapat99,Pattison00}

\section{A weighted form %representation
of hitting times for random walks}
\label{s_simp-walk}

Consider the hitting time quasi-metric in the case of {\it random walks\/} on connected positively weighted undirected\footnote{On hitting times for random walks on directed graphs, we refer to~\cite{BoleyRanjanZhang10TR}.} graphs $G,$ when the transition matrix is defined by~\eqref{e_W-T}. The class of such walks coincides with that of irreducible {reversible Markov chains\/}~\cite[Section\;3.2]{AldousFill14}.

In this case, $\ppi=(\pi_1\cdc\pi_n)$ is obviously proportional to $\bm1^{\xz\tra}\xy W.$ Furthermore, the hitting time quasi-metric $m$ is \cite[p.\:32]{DezaDeza11} a {\it weightable quasi-metric} (see also \cite{DezaDezaVidali12} and \cite[Chapter\:16]{DezaDeza16EShort}), i.e., there exists a {\it weight function} $u\!:\!V\to\mathbb{R}_{\ge 0}$
%$\bm u=(u_1, u_2\cdc u_n)$,
such that for all $i,j\in V$ it holds that

\beq{e_wei}%+
m(i,j)+u_i=
m(j,i)+u_j,
\eeq
where $u_i=u(i).$ %$,\; i\in V.$
Consequently, the hitting time quasi-metric $m$ has the {\em cyclic tour property}
(also called the {\em relaxed symmetry property}\/): for any $i, j, k\in V,$ it holds that
\beq{e_CTP}%+
m(i,j) + m(j,k) + m(k,i) = m(i,k) + m(k,j) + m(j,i).
\eeq

For unweighted graphs $G,$ this property appeared in~\cite[Lemma~2]{CoppersmithTetaliWinkler93}.
In turn, the cyclic tour property implies the weightability of~$m.$ Indeed, for an arbitrary $k\in V,$ set
\beq{e_weig}%+
u_i = m(k,i) - m(i,k),\quad i=1\cdc n. %i\in V.
\eeq
Now for any $i,j\in V,$ \eqref{e_CTP} gives
$
  m(i,j) - m(j,i) =
- m(j,k) - m(k,i) + m(i,k) + m(k,j) =
- u_i + u_j,
$
yielding \eqref{e_wei}. It remains to apply (if necessary) a shift that provides the defined function $u$ with non-negativity. Thus, \eqref{e_wei} and \eqref{e_CTP} are equivalent, and the weightability of hitting times for random walks on undirected weighted graphs follows from the cyclic tour property of any reversible Markov chain~\cite{AldousFill14}. Conversely, the cyclic tour property implies reversibility \cite{Tetali94} and thus, representability of the chain as a random walk on a weighted undirected graph. Hence, the weightability of hitting times indicates that the chain has the above representation.

A probabilistic interpretation of the weights $u_i$ is clear from \eqref{e_weig}: $u_1\cdc u_n$ are, up to a shift, {\em hitting time differences\/} from an arbitrary vertex $k$ to all vertices and back. They relatively measure {\em hitting asymmetry\/} of the vertices.
Corollary~\ref{c_hi0} supplies a structural description of this relative asymmetry: $u_i = f_{ki}/q_i - f_{ik}/q_k.$ It is worth recalling that the forest representation %of hitting times
involves the weighted digraph $\G$ of Theorem~\ref{thm_MFPT} rather than the initial graph~$G.$

The {\it commute time metric} $c$ on $V$ has now the representation
$$
c(i,j) = m(i,j) + m(j,i) = 2m(i, j) + u_i - u_j,
$$
while
$$
m(i,j) = \tfrac12\big(c(i,j) - u_i + u_j\big).
$$

In this case, the pair $(c, u)$ is a {\it weighted metric} on $V$, i.e., a metric with a weight function $u\!: V\to \mathbb{R}_{\ge 0}$ such that the {\it down-weighted condition} $c(i,j)\ge u_i-u_j$ is satisfied (\cite[Chapter\:6]{DezaDezaSikiric16}). %with a given weight $u_i$ on $i\in V$.
Furthermore, the function $p$,
$$
p(i, j) = m(i,j)+u_i = \tfrac12\big(c(i, j) + u_i + u_j\big),
$$
is a {\it partial metric} on $V$ (cf.\ \cite{{DezaDezaSikiric16}}), which means that for all $i, j, k\in V,$ it holds that:
\begin{enumerate}
\item $p(i, j)\ge 0;$
\item $p(i, j)\ge p(i, i)$ (\xz{\it small self-distances}\xy);
\item $p(i, i)=p(j, j)=p(i, j)$ $\Rightarrow$ $i=j$ (\xz{\it separation axiom}\xy);
\item $p(i, j)=p(j, i)$ (\xz{\it symmetry}\xy);
\item $p(i, j)\le p(i, k)+p(k, j)-p(k, k)$ (\xz{\it sharp triangle inequality}\xy).
\end{enumerate}

It is straightforward to check that
$$
\tfrac12\big(c(i, k) + c(k, j) - c(i, j)\big) =
             p(i, k) + p(k, j) - p(i, j) - p(k, k) =
             m(i, k) + m(k, j) - m(i, j),
$$
i.e., the respective triangle inequalities are equivalent on all three levels: of the weighted metric $c$, of the partial metric $p$, and of the weightable quasi-metric~$m$.

Moreover,
$$
m(i, j)\ge 0 \;\Leftrightarrow\; c(i, j) \ge u_i-u_j \;\Leftrightarrow\; p(i, j)\ge p(i, i).
$$
So, the {\it non-negativity condition} $m(i, j)\ge 0$ for the (weightable) quasi-metric $m$ is equivalent to the {\it down-weighted condition} $c(i, j) \ge u_i-u_j$ for the weighted metric $c$, and to the {\it small self-distances condition} $p(i, j)\ge p(i, i)$ for the partial metric $p$.

Now let us call a weightable quasi-metric $v$ along with weight function $u$ a {\it strong weighted quasi-metric} if for all $i, j\in V,$ $v(i, j)\le u_j$ holds.
Similarly, call a weighted metric $(d,u)$ a {\it strong weighted metric\/} if for all $i, j\in V,$ %$d(i, j)\ge u_i-u_j$ and
$d(i, j)\le u_i+u_j$ holds,
i.e., if it is not only {\it down-weighted}\/, but also {\it up-weighted}.
Finally, call a partial metric $p$ a {\it strong partial metric} if the {\it large self-distance condition} holds: $p(i, j)\le p(i, i)+p(j, j)$ for all $i, j\in V$.

It can be observed that
$$
m(i, j) \le u_j       \;\Leftrightarrow\;
c(i, j) \le u_i + u_j \;\Leftrightarrow\;
p(i, j) \le p(i, i) + p(j, j),\quad i,j\in V.
$$

Thus, the weightable quasi-metric $m$ with weight function $u$ is a {\it strong weighted quasi-metric} if and only if the weighted metric $(c, u)$ is a {\it strong weighted metric}, and if and only if the partial metric $p$ is a {\it strong partial metric}.

In this case, the strong weighted metric $(c, u)$ has an additional nice property.
Consider %the strong weighted metric $(c,u)$ on $V=\{1\cdc n\}$ as
the $(n+1)\!\times\!(n+1)$ matrix $[c'_{ij}]$, $0\le i,j\le n$, with $c'_{00}=0$,
$c'_{0i}=c'_{i0}=u_i$ for $i  \in V$, and
$c'_{ij}=c(i, j)$     for $i,j\in V$. In other words, the weight $u_i$ is considered as a distance from the point $i\in V$ to an additional point $0$: $u_i=c'(i, 0)=c'(0, i)$. In the case of \emph{strong\/} weighted metric $c$, the function $c'$ turns out to be a metric, since the addition of vertex $0$ does not violate the triangle inequality:
$c'(i, j)\le c'(i,0) + c'(0,j)$ and $c'(i, 0)\le c'(i,j) + c'(j,0).$

\medskip
The results presented in this paper demonstrate fruitful connections between the forest representation of hitting times and their metric properties.

\section{Examples}
\label{s_ex}

In this section, we illustrate the above concepts and results %, including Theorem~\ref{thm_MFPT}, Remark\;\ref{r_recu}, and various metric properties,
by two examples.

\subsection{Example 1: hitting times and their forest expression}
\label{ss_ex1}
Consider the Markov chain with transition matrix $T$ and the Laplacian matrix defined by\:\eqref{e_LaplP}: %$L$

$$
T=
\left[\begin{array}{rrrr}
      0&      1&      0&      0\\
      0&\frac45&\frac1{5\ms}& 0\\
\frac25&      0&\frac1{5\ms}&\frac25\\
      0&      0&\frac14&\frac34\\
\end{array}\right];\quad
%\end{bmatrix};\quad
L=
\left[\begin{array}{rrrr}
       1&     -1&            0&            0\\
       0&\frac15&-\frac1{5\ms}&            0\\
-\frac25&      0& \frac45     &-\frac2{5\ms}\\
       0&      0&-\frac14     & \frac14     \\
\end{array}\right].
$$

First, let us obtain the matrix $\M$ of hitting times by the direct use of~\eqref{Meyer3.3}.
Finding
$%\beq{exam_pi0}%-
\ppi=\tfrac1{25}(2, 10, 5, 8)
$ %\eeq
and
$$%   \beq{exam_L+0}%-
   L^\#=
\frac1{625}\left[\begin{array}{rrrr}
    463&  1065& -280& -1248\\
   -112&  1315& -155& -1048\\
    138&  -560&  470&   -48\\
    -62& -1560&  -30&  1652\\
   \end{array}\right]
$$%   \eeq
and substituting these in \eqref{Meyer3.3} yields
   \beq{M0}%+
   \M=
   \frac12\left[\begin{array}{rrrr}
   25 &  2 & 12 &   29\\ %@%
   23 &  5 & 10 &   27\\
   13 & 15 & 10 &   17\\
   21 & 23 &  8 & 6\frac14\\
   \end{array}\right].
   \eeq
Mention that $L^\#$ can be calculated (see, e.g., \cite[(i) of Proposition~15]{CheAga02ap}) by applying %& give a better reference!(?)
$$
L^\#=(L+\bm1\ppi)^{-1}-\bm1\ppi.
$$

\begin{figure}[htb] %[h!]
\begin{center}
% Drawing generated by LaTeX-CAD 1.9 - requires latexcad.sty 
% (c) 1998 John Leis leis@usq.edu.au 
\begin{picture}(25,28.5)
\thinlines
%\thicklines
%\drawcenteredtext{7.0}{25.0}{1}
%\drawcenteredtext{18.29}{25.0}{2}
%\drawcenteredtext{12.64}{17.48}{3}
%\drawcenteredtext{12.64}{6.44}{4}
\drawvector{9.19}{25.02}{7.19}{1}{0}
\drawvector{17.89}{22.73}{3.29}{-3}{-4}
\drawvector{10.69}{18.29}{3.4}{-3}{4}
\drawvector{13.1}{15.0}{6.3}{0}{-1}
\drawvector{12.0}{8.69}{6.3}{0}{1}
{\footnotesize
\drawcenteredtext{12.6}{26.7}{1}
% \drawcenteredtext{18.2}{19.89}{0.2}
% \drawcenteredtext{7.09}{19.89}{0.4}
% \drawcenteredtext{15.39}{11.69}{0.4}
% \drawcenteredtext{9.19}{11.69}{0.25}
\drawcenteredtext{7.14}{19.95}{$\frac25$}%{0.4}
\drawcenteredtext{18.15}{19.95}{$\frac15$}%{0.2}
\drawcenteredtext{9.85}{11.69}{$\frac14$}%{0.25}
\drawcenteredtext{15.31}{11.69}{$\frac25$}%{0.4}
}
\drawcircle{7.03}{25.04}{3.0}{1}
\drawcircle{18.28}{25.04}{3.0}{2}
\drawcircle{12.64}{17.52}{3.0}{3}
\drawcircle{12.64}{6.65}{3.0}{4}
\end{picture}
\end{center}\vspace{-3.2em} \caption{A weighted digraph corresponding to the Markov chain.\label{figG}}
\end{figure}
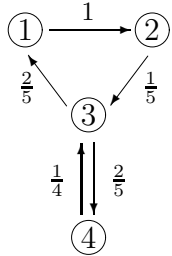

Now let us obtain $\M$ by means of Theorem~\ref{thm_MFPT}. The weighted digraph $\G$ without loops corresponding to the Markov chain under consideration is shown in Fig.~\ref{figG}. The converging trees of $\G$ are shown in Fig.~\ref{figT}, where the roots
are given in a boldface font.
\unitlength 1.50mm

%\vspace{-0.1em}
\unitlength 1.3mm
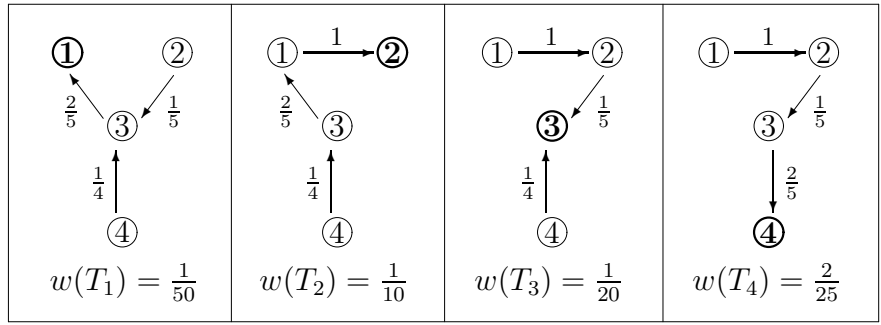
\begin{figure}[htb] %[h!]
\begin{center}
% Drawing generated by LaTeX-CAD 1.9 - requires latexcad.sty 
% (c) 1998 John Leis leis@usq.edu.au 
\begin{picture}(98,35.5) %100,37
\thinlines
\drawvector{22.87}{27.71}{3.27}{-3}{-4}
\drawvector{15.68}{23.29}{3.4}{-3}{4}
\drawvector{17.0}{13.68}{6.3}{0}{1}
{\footnotesize
\drawcenteredtext{12.39}{24.00}{$\frac25$}%{0.4}
\drawcenteredtext{22.90}{24.00}{$\frac15$}%{0.2}
\drawcenteredtext{15.18}{16.68}{$\frac14$}%{0.25}
}
\thicklines
\drawcircle{12.03}{30.04}{3.0}{\bf 1}
\thinlines
\drawcircle{23.28}{30.04}{3.0}{2}
\drawcircle{17.63}{22.52}{3.0}{3}
\drawcircle{17.63}{11.64}{3.0}{4}
\drawvector{36.15}{30.02}{7.19}{1}{0}
\drawvector{37.65}{23.29}{3.4}{-3}{4}
\drawvector{38.97}{13.68}{6.3}{0}{1}
{\footnotesize
\drawcenteredtext{39.58}{31.7}{1}
\drawcenteredtext{34.35}{24.00}{$\frac25$}%{0.4}
\drawcenteredtext{37.15}{16.68}{$\frac14$}%{0.25}
}
\drawcircle{34.01}{30.04}{3.0}{1}
\thicklines
\drawcircle{45.25}{30.04}{3.0}{\bf 2}
\thinlines
\drawcircle{39.61}{22.52}{3.0}{3}
\drawcircle{39.61}{11.63}{3.0}{4}
\drawvector{58.13}{30.02}{7.19}{1}{0}
\drawvector{66.83}{27.71}{3.27}{-3}{-4}
\drawvector{60.95}{13.68}{6.3}{0}{1}
{\footnotesize
\drawcenteredtext{61.54}{31.7}{1}
\drawcenteredtext{66.85}{24.00}{$\frac15$}%{0.2}
\drawcenteredtext{59.13}{16.68}{$\frac14$}%{0.25}
}
\drawcircle{55.97}{30.04}{3.0}{1}
\drawcircle{67.22}{30.04}{3.0}{2}
\thicklines
\drawcircle{61.59}{22.52}{3.0}{\bf 3}
\thinlines
\drawcircle{61.59}{11.63}{3.0}{4}
\drawvector{80.29}{30.02}{7.19}{1}{0}
\drawvector{88.98}{27.71}{3.27}{-3}{-4}
\drawvector{84.19}{20.0}{6.3}{0}{-1}
{\footnotesize
\drawcenteredtext{83.69}{31.7}{1}
\drawcenteredtext{89.00}{24.00}{$\frac15$}%{0.2}
\drawcenteredtext{86.10}{16.68}{$\frac25$}%{0.4}
}
\drawcircle{78.12}{30.04}{3.0}{1}
\drawcircle{89.37}{30.04}{3.0}{2}
\drawcircle{83.73}{22.52}{3.0}{3}
\thicklines
\drawcircle{83.73}{11.63}{3.0}{\bf 4}
\thinlines
\drawcenteredtext{17.82}{6.52}{$w(T_1)=\frac1{50}$}%0.02$}
\drawcenteredtext{39.27}{6.52}{$w(T_2)=\frac1{10}$}%0.1$}
\drawcenteredtext{61.25}{6.52}{$w(T_3)=\frac1{20}$}%0.05$}
\drawcenteredtext{83.56}{6.52}{$w(T_4)=\frac2{25}$}%0.08$}
\drawpath  {6.0}{35.0} {95.0}{35.0}
\drawpath {95.0}{35.0} {95.0} {2.5}
\drawpath {95.0} {2.5}  {6.0} {2.5}
\drawpath  {6.0} {2.5}  {6.0}{35.0}
\drawpath {28.8}{35.0} {28.8} {2.5}
\drawpath {50.6}{35.0} {50.6} {2.5}
\drawpath{72.94}{35.0}{72.94} {2.5}
\end{picture}
\end{center}\vspace{-2.5em} \caption{The converging trees $T_1,T_2,T_3$, and $T_4$ of $\G$.\label{figT}}
\end{figure}
\unitlength 1.50mm

Having the weights of these trees, by the definition of $q_i$ given in Section~\ref{s_forest} we obtain:
\beq{q_exam}%+
(q_1,q_2,q_3,q_4)
=\big(w(\{T_1\}),w(\{T_2\}),w(\{T_3\}),w(\{T_4\})\big)
=\tfrac1{100}(2,10,5,8).
\eeq
Since $q=\sum_{k=1}^4q_i=w(\{T_1,T_2,T_3,T_4\})=\si_3=\frac14$, \eqref{q_exam} implies
$%\beq{tilde_q_exam}%-
%\tilde q=(\tilde q_1,\tilde q_2,\tilde q_3,\tilde q_4)=
\frac{(q_1,q_2,q_3,q_4)}q%=(0.08,0.4,0.2,0.32)
                          =\tfrac1{25}(2,10,5,8).
$ %\eeq

In concordance with the Markov Chain Tree Theorem, this vector coincides with $\ppi$, the normalized left Perron vector of~$T.$

The 2-tree in-forests of $\G$ are shown in Fig.~\ref{figF}; the roots are given in a boldface font.

\unitlength 1.3mm
\begin{figure}[htb] %[h!]
\begin{center}
% Drawing generated by LaTeX-CAD 1.9 - requires latexcad.sty 
% (c) 1998 John Leis leis@usq.edu.au 
\begin{picture}(99,68) %,69
\thinlines
\drawcenteredtext{18.02}{38.04}{$w(F_1)=\frac2{25}$}%0.08$}
\drawvector{23.07}{59.58}{3.27}{-3}{-4}
\drawvector{15.89}{55.15}{3.4}{-3}{4}
{\footnotesize
\drawcenteredtext{12.61}{56.10}{$\frac25$}%{0.4}
\drawcenteredtext{23.13}{56.10}{$\frac15$}%{0.2}
}
\thicklines
\drawcircle{12.26}{61.9}{3.0}{\bf 1}
\drawcircle{17.86}{43.51}{3.0}{\bf 4}
\thinlines

\drawcircle{23.51}{61.9}{3.0}{2}
\drawcircle{17.86}{54.38}{3.0}{3}
\drawcenteredtext{40.16}{38.04}{$w(F_2)=\frac1{20}$}%0.05$}
\drawvector{45.22}{59.58}{3.27}{-3}{-4}
\drawvector{39.37}{45.54}{6.3}{0}{1}
{\footnotesize
\drawcenteredtext{45.28}{56.10}{$\frac15$}%{0.2}
\drawcenteredtext{37.55}{48.54}{$\frac14$}%{0.25}
}
\thicklines
\drawcircle{34.4}{61.9}{3.0}{\bf 1}
\drawcircle{40.01}{54.38}{3.0}{\bf 3}
\thinlines

\drawcircle{45.65}{61.9}{3.0}{2}
\drawcircle{40.01}{43.51}{3.0}{4}
\drawcenteredtext{61.97}{38.04}{$w(F_3)=\frac1{10}$}%0.1$}
\drawvector{59.84}{55.15}{3.4}{-3}{4}
\drawvector{61.16}{45.54}{6.3}{0}{1}
{\footnotesize
\drawcenteredtext{56.56}{56.10}{$\frac25$}%{0.4}
\drawcenteredtext{59.36}{48.54}{$\frac14$}%{0.25}
}
\thicklines
\drawcircle{56.2}{61.9}{3.0}{\bf 1}
\drawcircle{67.45}{61.9}{3.0}{\bf 2}
\thinlines
\drawcircle{61.81}{54.38}{3.0}{3}
\drawcircle{61.81}{43.51}{3.0}{4}
\drawcenteredtext{83.48}{38.63}{$w(F_4)=\frac25$}%0.4$}
\drawvector{80.37}{61.96}{7.19}{1}{0}
\drawvector{81.87}{55.23}{3.4}{-3}{4}
{\footnotesize
\drawcenteredtext{83.8}{63.63}{1}
\drawcenteredtext{78.57}{56.10}{$\frac25$}%{0.4}
}
\drawcircle{78.22}{61.98}{3.0}{1}
\drawcircle{83.83}{54.46}{3.0}{3}
\thicklines
\drawcircle{89.47}{61.98}{3.0}{\bf 2}
\drawcircle{83.83}{43.58}{3.0}{\bf 4}
\thinlines
\drawcenteredtext{17.57}{6.58}{$w(F_5)=\frac14$}%0.25$}
\drawvector{14.46}{29.9}{7.19}{1}{0}
\drawvector{17.27}{13.56}{6.3}{0}{1}
{\footnotesize
\drawcenteredtext{17.88}{31.58}{1}
\drawcenteredtext{15.46}{16.56}{$\frac14$}%{0.25}
}
\drawcircle{12.32}{29.91}{3.0}{1}
\thicklines
\drawcircle{23.56}{29.91}{3.0}{\bf 2}
\drawcircle{17.91}{22.4}{3.0}{\bf 3}
\thinlines
\drawcircle{17.91}{11.52}{3.0}{4}
\drawcenteredtext{39.83}{6.4}{$w(F_6)=\frac15$}%0.2$}
\drawvector{36.56}{29.9}{7.19}{1}{0}
\drawvector{45.26}{27.61}{3.27}{-3}{-4}
{\footnotesize
\drawcenteredtext{39.97}{31.58}{1}
\drawcenteredtext{45.58}{24.76}{$\frac15$}%{0.2}
}
\drawcircle{34.4}{29.93}{3.0}{1}
\drawcircle{45.65}{29.93}{3.0}{2}
\thicklines
\drawcircle{40.01}{22.4}{3.0}{\bf 3}
\drawcircle{40.01}{11.52}{3.0}{\bf 4}
\thinlines
\drawcenteredtext{61.53}{6.4}{$w(F_7)=\frac2{25}$}%0.08$}
\drawvector{66.96}{27.61}{3.27}{-3}{-4}
\drawvector{62.17}{19.88}{6.3}{0}{-1}
{\footnotesize
\drawcenteredtext{67.27}{24.76}{$\frac15$}%{0.2}
\drawcenteredtext{64.16}{16.56}{$\frac25$}%{0.4}
}
\thicklines
\drawcircle{56.1}{29.93}{3.0}{\bf 1}
\drawcircle{61.71}{11.52}{3.0}{\bf 4}
\thinlines
\drawcircle{67.35}{29.93}{3.0}{2}
\drawcircle{61.71}{22.4}{3.0}{3}

\drawcenteredtext{83.5}{6.4}{$w(F_8)=\frac25$}%0.4$}
\drawvector{80.23}{29.9}{7.19}{1}{0}
\drawvector{84.13}{19.88}{6.3}{0}{-1}
{\footnotesize
\drawcenteredtext{83.63}{31.58}{1}
\drawcenteredtext{86.13}{16.56}{$\frac25$}%{0.4}
}
\drawcircle{78.06}{29.93}{3.0}{1}
\drawcircle{83.68}{22.4}{3.0}{3}
\thicklines
\drawcircle{89.31}{29.93}{3.0}{\bf 2}
\drawcircle{83.68}{11.52}{3.0}{\bf 4}
\thinlines
\drawpath {7.0}{68.0}{7.0} {4.0}
\drawpath {7.0}{4.0} {95.0}{4.0}
\drawpath{95.0}{4.0} {95.0}{68.0}
\drawpath{95.0}{68.0}{7.0} {68.0}
\drawpath{7.0} {34.0}{95.0}{34.0}
\drawpath{29.0}{68.0}{29.0}{4.0}
\drawpath{51.0}{68.0}{51.0}{4.0}
\drawpath{73.0}{68.0}{73.0}{4.0}
\end{picture}
\end{center}
\vspace{-3.0em}\caption{The 2-tree in-forests $F_1\cdc F_8$ of $\G$.\label{figF}}
\end{figure}
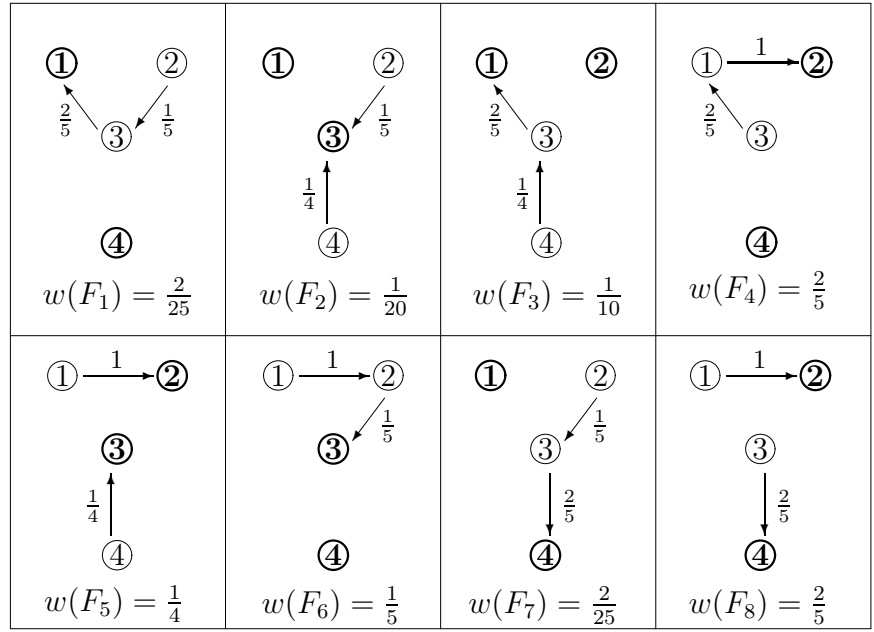
\unitlength 1.50mm

In Theorem~\ref{thm_MFPT}, $f_{ij}$ is defined as the total weight of 2-tree in-forests of $\G$ that have one tree containing $i$ and the other tree converging to~$j$. Therefore,
   \beqq\nonumber
   [f_{ij}]
   &=&
   \left[\begin{array}{llll}
   0                   &w(\{F_3\})             &w(\{F_2,F_5\}) &w(\{F_1,F_4,F_6,F_7,F_8\})\\
   w(\{F_2,F_3,F_7\})  &0                      &w(\{F_5\})     &w(\{F_1,F_4,F_6,F_8\})    \\
   w(\{F_2,F_7\})      &w(\{F_3,F_5,F_8\})     &0              &w(\{F_1,F_4,F_6\})        \\
   w(\{F_1,F_2,F_7\})  &w(\{F_3,F_4,F_5,F_8\}) &w(\{F_6\})     &0                         \\
   \end{array}\right]
   \\
   &=&
   \frac1{100}\left[\begin{array}{rrrr}
    0 & 10 &30 &116\\ %@%
   23 &  0 &25 &108\\
   13 & 75 & 0 & 68\\
   21 &115 &20 &  0\\
   \end{array}\right],
   \label{exam_fij}%+
   \eeqq
where $w(A)$ is the weight of a set $A$ of digraphs. Moreover, $\si_2=w(\{F_1\cdc F_8\})=\frac{39}{25}.$

Substituting \eqref{q_exam}--\eqref{exam_fij} in \eqref{MFPT_forest} yields the matrix $\M$ of hitting times coinciding with \eqref{M0}.

Remark\;\ref{r_recu} enables one to avoid generating the converging trees and 2-tree in-forests of~$\G$.
%shown in Fig.~\ref{figT} and Fig.~\ref{figF}.
Instead, $f_{ij}$ and $q_j$ can be computed by means of the recurrent procedure \eqref{req1}--\eqref{req2}. Starting with $Q_0=I$, $\sigma_0=1,$ for this example we have:

\Up{1.6}\beqq
%   \si_1&=&\frac{\tr(LQ_0)}{1}=2.25,\nonumber\\
   Q_1  &=&-LQ_0+\frac{\tr(LQ_0)}1I=
                                 \frac1{20}\left[\begin{array}{rrrr}
                                   25&   20&    0&   0\\
                                    0&   41&    4&   0\\
                                    8&    0&   29&   8\\
                                    0&    0&    5&  40\\
                                 \end{array}\right],\quad \si_1=\frac94\;;\nonumber\\
%-%   \eeqq
%-%\Up{1.5}\beqq
%   \si_2&=&\frac{\tr(LQ_1)}{2}=1.56,\nonumber\\
   Q_2  &=&-LQ_1+\frac{\tr(LQ_1)}2I=
                                 \frac1{100}\left[\begin{array}{rrrr}
                                 31& 105& 20&   0\\
                                  8& 115& 25&   8\\
                                 18&  40& 50&  48\\
                                 10&   0& 30& 116\\
                                 \end{array}\right],\quad \si_2=\frac{39}{25}\;;\label{examQ2}\\%+
%-%   \eeqq
%-%\Up{1.5}\beqq
%   \si_3&=&\frac{\tr(LQ_2)}{3}=0.25,\nonumber\\
   Q_3  &=&-LQ_2+\frac{\tr(LQ_2)}3I
         =\frac1{100}\,\bm1
%         \left[\begin{array}{rrrr}2& 10& 5& 8\end{array}\right]
          \begin{bmatrix}\,2& \!10& \!5& \!8\,\end{bmatrix},\quad \si_3=\frac14\;.\label{examQ3}%+
   \eeqq

%Using $Q_2=[q^{(2)}_{ij}]_{4\times 4}$
By \eqref{fijdef} we have $f_{ij}=q^{(2)}_{jj}-q^{(2)}_{ij}$, $i,j=1\cdc 4.$
Thereby \eqref{examQ2} provides the matrix $[f_{ij}]_{4\times 4},$ which coincides with~\eqref{exam_fij}.
Eq.~\eqref{examQ3} yields $(q_1,q_2,q_3,q_4)=\frac1{100}(2,10,5,8)$, which coincides with~\eqref{q_exam}. Now using Theorem~\ref{thm_MFPT} we obtain the matrix~\eqref{M0} of hitting times again.

By \eqref{examQ2}, \eqref{examQ3}, and Corollary~\ref{c_kemeny}, Kemeny's constant of this chain is $K=1+\frac{\si_2}{\si_3}=7\frac6{25}$ and by Remark\;\ref{r_mvsm}, $\sum_{j=1}^4\pi_j\xy m_{ij}=\frac{\si_2}{\si_3}=6\frac6{25}$ for any $i=1\cdc 4.$

It is easy\x{} to observe that the hitting time quasi-metric $m$ defined by the matrix
\Up{0.5}
$$% \beq{M0A}%-
   M=
   \frac12\left[\begin{array}{rrrr}
    0 &   2 & 12 & 29\\ %@%
   23 &   0 & 10 & 27\\
   13 &  15 &  0 & 17\\
   21 &  23 &  8 &  0\\
   \end{array}\right]
$$%   \eeq
is cutpoint additive.
For example,
$m(4,3) + m(3,2) = m(4,2),$
$m(1,3) + m(3,4) = m(1,4),$
$m(2,3) + m(3,4) = m(2,4),$ however,
$m(3,2) + m(2,4) > m(3,4).$

On the other hand, it is not weightable, as the cyclic tour property is violated:

\Up{1.0}
$$
m(1,2) + m(2,3) + m(3,1)\neq m(1,3) + m(3,2) + m(2,1).
$$

The corresponding commute time metric $c$, $c(i,j)=m(i,j)+m(j,i)$ is defined by %the matrix
$$%\beq{M11}%-
C=\frac12\left[
\begin{array}{rrrr}
    0 & 25 & 25 & 49\\ %@%
   25 &  0 & 25 & 49\\
   25 & 25 &  0 & 25\\
   49 & 49 & 25 &  0\\
   \end{array}
\right].
$$%\eeq

\subsection{Example 2: hitting metric functions for an undirected graph}
\label{ss_ex3}
To illustrate the concept of weighted metric on a nontrivial example, consider a %simple
random walk on the undirected unweighted graph $G=(V,E)$ with $V=\{1\cdc6\}$ and $E=\{\{1,2\},\{2,3\},\{3,4\},$ $\{3,5\},$ $\{4,5\},\{5,6\}\}$, whose automorphism group is trivial (Fig.\;4).
\begin{figure}[htb] %[h!]
\begin{center}
%\Up{3}
% Drawing generated by LaTeX-CAD 1.9 - requires latexcad.sty 
% (c) 1998 John Leis leis@usq.edu.au 
\begin{picture}(45,17)
\thinlines
%\thicklines
\drawpath{ 8.50}{6.00}{12.50}{ 6.00}
\drawpath{15.50}{6.00}{19.50}{ 6.00}
\drawpath{22.50}{6.00}{26.50}{ 6.00}
\drawpath{29.50}{6.00}{33.50}{ 6.00}
\drawpath{21.00}{7.50}{21.00}{11.50}
\drawpath{22.13}{11.87}{26.88}{ 7.12}
{\footnotesize }
\drawcircle{ 7.00}{ 6.00}{3.0}{1}
\drawcircle{14.00}{ 6.00}{3.0}{2}
\drawcircle{21.00}{ 6.00}{3.0}{3}
\drawcircle{21.00}{13.00}{3.0}{4}
\drawcircle{28.00}{ 6.00}{3.0}{5}
\drawcircle{35.00}{ 6.00}{3.0}{6}
\end{picture}
\end{center}\vspace{-2.6em} \caption{A graph $G$ whose automorphism group is trivial.\label{figG13}}
\end{figure}
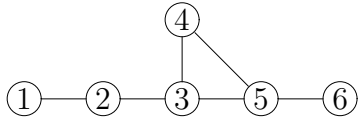
\unitlength 1.50mm

Define the transition matrix of the corresponding Markov chain by \eqref{e_W-T}.
As the vertex degrees are $d(1)=d(6)=1$, $d(2)=d(4)=2$, and $d(3)=d(5)=3,$ it holds that

$$
T=
\left[\begin{array}{rrrrrr}
  0&   1&   0&   0&   0&   0\\
1/2&   0& 1/2&   0&   0&   0\\
  0& 1/3&   0& 1/3& 1/3&   0\\
  0&   0& 1/2&   0& 1/2&   0\\
  0&   0& 1/3& 1/3&   0& 1/3\\
  0&   0&   0&   0&   1&   0
\end{array}\right].
$$

$T$ is the matrix of arc weights of the corresponding digraph $\G$ without loops. %with $V=\{1\cdc 6\}$.
%The corresponding digraph $\G$ without loops on the same vertex set has arc weights $w(1,2)=1$, $w(2,k)=1/2$, $k=1, 3$, $w(3,k)=1/3$, $k=2, 4, 5$, $w(4,k)=1/2$, $k=3, 5$, $w(5,k)=1/3$, $k=3, 4, 6$, and $w(6,5)=1$.
%Now we compute the matrix $\M$ of hitting times by means of Theorem~\ref{thm_MFPT}.
As well as for a general Markov chain, the matrix $\M$ of hitting times can be computed using Theorem~\ref{thm_MFPT}. However, in the present case, there is no need to enumerate trees for obtaining
\beq{q_exam13}%+
(q_1,q_2,q_3,q_4, q_5, q_6)=\tfrac1{12}\big(1,2,3,2,3,1\big),
\eeq
$q=\sum_{k=1}^6q_i=\si_5=1$ and finally $\ppi=q^{-1}{(q_1\cdc q_6)}=\tfrac1{12}(1,2,3,2,3,1),$ as we know that for such %this kind of
random walks, $\ppi$ is proportional to $\bm1^{\xz\tra}\xy W,$ where $W$ is the edge weight matrix of~$G.$ %the initial graph~$G.$
%Enumerating 18 converging trees of $\G$ and using the definition of $q_i$ one finds the row vector whose entries are proportional to the degrees of the corresponding vertices.
%Since $q=\sum_{k=1}^6q_i=1$, \eqref{q_exam13} implies $q^{-1}{(q_1,q_2,q_3,q_4,q_5,q_6)}=(q_1,q_2,q_3,q_4,q_5,q_6).$
%This vector coincides with $\ppi$, the normalized left Perron vector of~$T$.

Using all 76 2-tree in-forests of $\G$ one may obtain the matrix $[f_{ij}],$ where $f_{ij}$ is the total weight of 2-tree in-forests of $\G,$ where one tree contains $i$ and the other converges to~$j$: %In Theorem~\ref{thm_MFPT},

\Up{0.6}
\beq{exam_fij1A3}%+
   [f_{ij}]=
   \frac{1}{36}
   \left[\begin{array}{rrrrrr}
    0&  6& 36& 56& 78& 59\\ %@%
   33&  0& 27& 50& 69& 56\\
   60& 54&  0& 32& 42& 47\\
   68& 70& 24&  0& 30& 43\\
   70& 74& 30& 28&  0& 33\\
   73& 80& 39& 34&  9&  0\\
   \end{array}\right];\quad \si_4=\sum_{j=1}^6 f_{ij}=\frac{235}{36},\:\;i=1\cdc6.
   \eeq

Kemeny's constant is $1+\frac{\si_4}{\si_5}\xz=\xz7\frac{19}{36}.$ Substituting \eqref{q_exam13}--\eqref{exam_fij1A3} into \eqref{MFPT_forest} and \eqref{MFPT_forestA} yields the matrix $\M$ of hitting times and the cutpoint additive quasi-metric $m$ represented by matrix~$M$:

\Up{0.8}
$$%\beq{M113}%+
\M= \frac13
\left[\begin{array}{rrrrrr}
36&  3& 12& 28& 26& 59\\ %@%
33& 18&  9& 25& 23& 56\\
60& 27& 12& 16& 14& 47\\
68& 35&  8& 18& 10& 43\\
70& 37& 10& 14& 12& 33\\
73& 40& 13& 17&  3& 36\\
\end{array}\right];\quad
%$$%\beq{oldM11}
M=\frac13
\left[
\begin{array}{rrrrrr}
 0&  3& 12& 28& 26& 59\\ %@%
33&  0&  9& 25& 23& 56\\
60& 27&  0& 16& 14& 47\\
68& 35&  8&  0& 10& 43\\
70& 37& 10& 14&  0& 33\\
73& 40& 13& 17&  3&  0\\
\end{array}
\right].
$$%\eeq

Furthermore, $m$ is a weightable quasi-metric, whose (non-negative and defined up to a positive shift) weight function is defined by the row vector
$\bm u=\tfrac13(48, 18, 0, 8, 4, 34).$

The corresponding commute time metric $c$, $c_{ij}=m(i, j)+m(j, i),$
the resistance distance $\Omega(i,j)=(\bm1^\tra W\bm1)^{-1}c(i,j),$
and the partial metric $p(i, j)=\frac{c(i, j)+u_i+u_j}{2}$, are given by %the matrices:

\Up{1.3}
$$%\beq{oldM11}%-
C=\left[
\begin{array}{rrrrrr}
 0& 12& 24& 32& 32& 44\\ %@%
12&  0& 12& 20& 20& 32\\
24& 12&  0&  8&  8& 20\\
32& 20&  8&  0&  8& 20\\
32& 20&  8&  8&  0& 12\\
44& 32& 20& 20& 12&  0\\
\end{array}
\right];\;\;\Omega=\tfrac1{12}C;\;\; %[\Omega_{ij}]
%$$%\eeq
%$$%\beq{oldM11}
P=\frac{1}{3}\left[
\begin{array}{rrrrrr}
 48& 51& 60& 76& 74& 107\\ %@%
 51& 18& 27& 43& 41&  74\\
 60& 27&  0& 16& 14&  47\\
 76& 43& 16&  8& 18&  51\\
 74& 41& 14& 18&  4&  37\\
107& 74& 47& 51& 37&  34\\
\end{array}
\right].
$$%\eeq

\Up{0.6}
For any weight function $u$ such that $m(i, j)\le u_j$, $i, j=1\cdc6$, i.e., starting from $\bm u=\frac{1}{3}(73,43,25,33,29,59)$, functions $m$, $c,$ and $p$ are strong on the corresponding level: $m$ is a strong weighted quasi-metric, $c$ a strong weighted metric, and $p$ a strong partial metric.

Moreover, in this case, the function $c'\!: V'\times V'\rightarrow \mathbb{R}$, where $V'\xz=\xz\{0, 1\cdc6\}$, $c'(0, 0)\xz=\xz0$, $c'(0, i)=c'(i, 0)=u_i$, and $c'(i, j)=c(i, j)$ for $i, j=1\cdc6$, is a metric on $V'.$ For $\bm u=\frac{1}{3}(73, 43, 25,  33, 29, 59)$ its matrix is:

\Up{0.6}
$$%\beq{M11b}%-
C'=\frac{1}{3}
\left[
\begin{array}{rrrrrrr}
 0&  73& 43& 25& 33& 29&  59\\
73&   0& 36& 72& 96& 96& 132\\%@%
43&  36&  0& 36& 60& 60&  96\\
25&  72& 36&  0& 24& 24&  60\\
33&  96& 60& 24&  0& 24&  60\\
29&  96& 60& 24& 24&  0&  36\\
59& 132& 96& 60& 60& 36&   0\\
\end{array}
\right].
$$%\eeq

%\nocite{FoussPirotte07simi}
%\nocite{Biggs97,Tetali91}
%\nocite{Bapat11Passage,CarmonaEncinasMitjana17}
%\nocite{BoleyRanjanZhang10TR}
%\nocite{CheSha95b}
\nocite{Sonin99}

\small
\bibliographystyle{endm}%{unsrt}%{sweunsrt}%{dk-unsrt}%{jurunsrt}%{namunsrt}%{is-unsrt}%{unsrt}%{amsplain}%{abbrv}%{plain}%{abbrvnat}
\bibliography{all2}

\begin{thebibliography}{10}
\expandafter\ifx\csname url\endcsname\relax
  \def\url#1{\texttt{#1}}\fi
\expandafter\ifx\csname urlprefix\endcsname\relax\def\urlprefix{URL }\fi
\newcommand{\enquote}[1]{``#1''}

\bibitem{AldousFill14}
Aldous, D. and J.~A. Fill, \enquote{Reversible {Markov} Chains and Random Walks
  on Graphs,} 2002, unfinished monograph,
  \url{http://www.stat.berkeley.edu/\textasciitilde aldous/RWG/book.html}.

\bibitem{BapatSivasubramanian11}
Bapat, R.~B. and S.~Sivasubramanian, \emph{{Identities for minors of the
  Laplacian, resistance and distance matrices}}, Linear Algebra and Its
  Applications \textbf{435} (2011), pp.~1479--1489.

\bibitem{BoleyRanjanZhang10TR}
Boley, D., G.~Ranjan and Z.-l. Zhang, \emph{Commute times for a directed graph
  using an asymmetric {Laplacian}}, Linear Algebra and its Applications
  \textbf{435} (2011), pp.~224--242.

\bibitem{Catoni99}
Catoni, O., \emph{Simulated annealing algorithms and {M}arkov chains with rare
  transitions}, in: \emph{S\'eminaire de Probabilit\'es, XXXIII},  LNM
  \textbf{1709}, Springer, Berlin, 1999 pp. 69--119.

\bibitem{CatralNeumannXu05}
Catral, M., M.~Neumann and J.~Xu, \emph{Proximity in group inverses of
  {M}-matrices and inverses of diagonally dominant {M}-matrices}, Linear
  Algebra and its Applications \textbf{409} (2005), pp.~32--50.

\bibitem{ChandraRaghavanRuzzoSmolenskyTiwari89}
Chandra, A.~K., P.~Raghavan, W.~L. Ruzzo, R.~Smolensky and P.~Tiwari, \emph{The
  electrical resistance of a graph captures its commute and cover times}, in:
  \emph{Proc. 21st Annual ACM Symp. on Theory of Computing} (1989), pp.
  574--586.

\bibitem{Che07mfpt}
Chebotarev, P., \emph{A graph theoretic interpretation of the mean first
  passage times}, arXiv preprint math.PR/0701359 (2007).

\bibitem{Che13Paris}
Chebotarev, P., \emph{Studying new classes of graph metrics}, in: F.~Nielsen
  and F.~Barbaresco, editors, \emph{Proceedings of the SEE Conference
  ``Geometric Science of Information'' (GSI-2013)}, Lecture Notes in Computer
  Science, LNCS 8085 (2013), pp. 207--214.

\bibitem{CheAga02ap}
Chebotarev, P. and R.~Agaev, \emph{{Forest matrices around the Laplacian
  matrix}}, Linear Algebra and its Applications \textbf{356} (2002),
  pp.~253--274.

\bibitem{CoppersmithTetaliWinkler93}
Coppersmith, D., P.~Tetali and P.~Winkler, \emph{Collisions among random walks
  on a graph}, SIAM Journal on Discrete Mathematics \textbf{6} (1993),
  pp.~363--374.

\bibitem{DezaDezaSikiric16}
Deza, E., M.~Deza and M.~D. Sikiri{\v{c}}, \enquote{Generalizations of Finite
  Metrics and Cuts,} World Scientific, 2016.

\bibitem{DezaDeza11}
Deza, M. and E.~Deza, \emph{Cones of partial metrics}, Contributions to
  Discrete Mathematics \textbf{6} (2011), pp.~26--47.

\bibitem{DezaDezaVidali12}
Deza, M., E.~Deza and J.~Vidali, \emph{Cones of weighted and partial metrics},
  in: \emph{Proceedings of the Internat. Conference on Algebra 2010: Advances
  in Algebraic Structures} (2012), pp. 177--197.

\bibitem{DezaDeza16EShort}
Deza, M.~M. and E.~Deza, \enquote{Encyclopedia of Distances,} Springer,
  Berlin--Heidelberg, 2016.

\bibitem{DoyleSnell84}
Doyle, P.~G. and J.~L. Snell, \enquote{Random Walks and Electric Networks,}
  Mathematical Association of America, Washington D. C., 1984.

\bibitem{EllensSpieksma11}
Ellens, W., F.~M. Spieksma, P.~Van~Mieghem, A.~Jamakovic and R.~E. Kooij,
  \emph{Effective graph resistance}, Linear Algebra and Its Applications
  \textbf{435} (2011), pp.~2491--2506.

\bibitem{FreidlinWentzell84}
Freidlin, M.~I. and A.~D. Wentzell, \enquote{Random Perturbations of Dynamical
  Systems,} Springer, New York, 1984.

\bibitem{GvishianiGurvich87En}
Gvishiani, A.~D. and V.~A. Gurvich, \emph{Metric and ultrametric spaces of
  resistances}, Russian Mathematical Surveys \textbf{42} (1987), pp.~235--236.

\bibitem{Hausdorff1927}
Hausdorff, F., \enquote{Grundz\"uge der Mengenlehre,} Walter de Gruyter,
  Berlin, 1927.

\bibitem{Hunter14}
Hunter, J.~J., \emph{The role of {Kemeny}'s constant in properties of {Markov}
  chains}, Communication in Statistics -- Theory and Methods \textbf{43}
  (2014), pp.~1309--1321.

\bibitem{Hunter16}
Hunter, J.~J., \emph{Accurate calculations of stationary distributions and mean
  first passage times in {Markov} renewal processes and {Markov} chains},
  Special Matrices \textbf{4} (2016), pp.~151--175.

\bibitem{KemenySnellKnapp76}
Kemeny, J.~G., J.~L. Snell and A.~W. Knapp, \enquote{Denumerable Markov Chains,
  volume 40 of Graduate Texts in Mathematics,} Springer-Verlag, New York, 1976.

\bibitem{KirklandZeng16}
Kirkland, S. and Z.~Zeng, \emph{Kemeny's constant and an analogue of {Braess'}
  paradox for trees}, Electronic Journal of Linear Algebra \textbf{31} (2016),
  pp.~444--464.

\bibitem{KirklandNeumann12book}
Kirkland, S.~J. and M.~Neumann, \enquote{{Group Iverses of M-matrices and Their
  Applications},} CRC Press, 2012.

\bibitem{KleinZhu98}
Klein, D. and H.~Zhu, \emph{Distances and volumina for graphs}, Journal of
  Mathematical Chemistry \textbf{23} (1998), pp.~179--195.

\bibitem{KleinRandic93}
Klein, D.~J. and M.~Randi\'c, \emph{Resistance distance}, Journal of
  Mathematical Chemistry \textbf{12} (1993), pp.~81--95.

\bibitem{LeightonRivest83}
Leighton, T. and R.~L. Rivest, \emph{The {Markov} chain tree theorem}, Computer
  Science Technical Report MIT/LCS/TM-249, Laboratory of Computer Science, MIT,
  Cambridge, Mass. (1983).

\bibitem{LeightonRivest86}
Leighton, T. and R.~L. Rivest, \emph{Estimating a probability using finite
  memory}, IEEE Transactions on Information Theory \textbf{32} (1986),
  pp.~733--742.

\bibitem{Meyer75}
Meyer{,}~Jr., C.~D., \emph{The role of the group generalized inverse in the
  theory of finite {Markov} chains}, SIAM Review \textbf{17} (1975),
  pp.~443--464.

\bibitem{OlivieriScoppola96}
Olivieri, E. and E.~Scoppola, \emph{{Markov chains with exponentially small
  transition probabilities: first exit problem from a general domain. II. The
  general case}}, Journal of Statistical Physics \textbf{84} (1996),
  pp.~987--1041.

\bibitem{PitmanTang16}
Pitman, J. and W.~Tang, \emph{{Tree formulas, mean first passage times and
  Kemeny's constant of a Markov chain}}, Bernoulli \textbf{24} (2018),
  pp.~1942--1972.

\bibitem{RaoMitra71}
Rao, C.~R. and S.~K. Mitra, \enquote{Generalized Inverse of Matrices and its
  Applications,} Wiley, New York, 1971.

\bibitem{SeshuReed61}
Seshu, S. and M.~B. Reed, \enquote{Linear Graphs and Electrical Networks,}
  Addison-Wesley, Reading, MA, 1961.

\bibitem{Shapiro87MathMag}
Shapiro, L.~W., \emph{{An electrical lemma}}, Mathematics Magazine \textbf{60}
  (1987), pp.~36--38.

\bibitem{Sharpe67a}
Sharpe, G.~E., \emph{Solution of the $(m+1)$-terminal resistive network problem
  by means of metric geometry}, in: \emph{Proceedings of the First Asilomar
  Conference on Circuits and Systems}, Pacific Grove, CA, 1967, pp. 319--328.

\bibitem{SharpeStyan67}
Sharpe, G.~E. and G.~P.~H. Styan, \emph{A note on equicofactor matrices},
  Proceedings of the IEEE \textbf{55} (1967), pp.~1226--1227.

\bibitem{Sonin99}
Sonin, I., \emph{The state reduction and related algorithms and their
  applications to the study of {Markov} chains, graph theory, and the optimal
  stopping problem}, Advances in Mathematics \textbf{145} (1999), pp.~159--188.

\bibitem{Tetali94}
Tetali, P., \emph{An extension of {Foster's} network theorem}, Combinatorics,
  Probability and Computing \textbf{3} (1994), pp.~421--427.

\bibitem{WentzellFreidlin70a}
Wentzell, A.~D. and M.~I. Freidlin, \emph{On small random perturbations of
  dynamical systems}, Russian Mathematical Surveys \textbf{25} (1970),
  pp.~1--55.

\bibitem{WentzellFreidlin79e}
Wentzell, A.~D. and M.~I. Freidlin, \enquote{Fluctuations in Dynamical Systems
  under Small Random Perturbations,} Nauka, Moscow, 1979, in Russian.

\bibitem{Wilson1931}
Wilson, W., \emph{On quasi-metric spaces}, American Journal of Mathematics
  \textbf{53} (1931), pp.~675--684.

\bibitem{YoungScardoviLeonard16}
Young, G.~F., L.~Scardovi and N.~E. Leonard, \emph{A new notion of effective
  resistance for directed graphs --- {Part I: Definition}s and properties},
  IEEE Transactions on Automatic Control \textbf{61} (2016), pp.~1727--1736.

\end{thebibliography}

\end{document}